# From Cell Decomposition to Motivic Integration, Hensel Minimality, and Point Counting

Raf Cluckers[†]

**Abstract.** Following Macintyre's quantifier elimination result for semi-algebraic $p$-adic sets from 1976, Denef's cell decomposition from 1984 has played a catalyzing role for both motivic integration and Hensel minimality. Motivic integration has in turn been applied in the Langlands program via transfer principles to change the characteristic of the local field, e.g., for the fundamental lemma. Hensel minimality has been used to study rational points on definable sets, providing non-archimedean analogues of results by Pila and Wilkie in o-minimal structures. I will review some related results and open questions.

**1 Introduction.** Soon after Macintyre [87] obtained his renowned quantifier elimination result for semi-algebraic $p$-adic sets in 1976, Denef [53] used it in 1984 to show cell decomposition for the semi-algebraic $p$-adic sets, which led to new rationality proofs of Igusa's local zeta function without relying on Hironaka's embedded resolution result. Indeed, $p$-adic cell decomposition served as an alternative to Hironaka's embedded resolution of singularities and enabled Denef to answer Serre's question [103] on rationality for a variant of Igusa's zeta function in [53].

Macintyre's result can be seen as the analogue for $\mathbb{Q}_p$ of Tarski's result [107] about real semi-algebraic sets. Write $P_m^K$ for the set of nonzero $m$th powers in a field $K$ with $m>0$ an integer. Call a set $X\subset K^n$ *semi-algebraic* if it can be obtained in finitely many steps from taking zero loci (in $K^\ell$ for some $\ell$) of polynomials over $K$, Boolean combinations, and images under coordinate projections $K^{m'+n'}\to K^{m'}$ with nonnegative integers $n,\ell,m',n'$. Now the results by Tarski and Macintyre can be summarized together as follows.

Theorem 1.1 (Tarski / Macintyre quantifier elimination [107], [87]). *Let $K$ be either $\mathbb{R}$ or $\mathbb{Q}_p$. Then a set $X\subset K^n$ is semi-algebraic if and only if it can be written as a Boolean combination of sets of the form*

$$\{x\in K^n \mid f(x)\in P_m^K\}$$

*for some polynomials $f$ in $K[X_1,\dots,X_n]$ and some integers $m>0$.*

Of course, for the reals only the squares are needed (or an order relation), but Theorem 1.1 clearly shows the similarities between the real and the $p$-adic cases. This quantifier elimination result enabled the study of real and $p$-adic semi-algebraic geometry, as it puts the focus on Boolean combinations of basic conditions and avoids the projections.

In order to understand the geometry of semi-algebraic sets, it is important to study the behavior of semi-algebraic functions on them. This is done via results on cell decomposition and preparation that we recall here in a basic, one-dimensional case. Recall that a function $f:X\to Y$ between semi-algebraic sets is called semi-algebraic if its graph is a semi-algebraic set.

[†]Univ. Lille, CNRS, UMR 8524 - Laboratoire Paul Painlevé, F-59000 Lille, France, and KU Leuven, Department of Mathematics, B-3001 Leuven, Belgium (Raf.Cluckers@univ-lille.fr, https://rcluckers.perso.math.cnrs.fr/).

Let $N>0$ be an integer. Note that $1+Np\mathbb{Z}_p$ is a subgroup of $\mathbb{Q}_p^\times$, and write $\sigma_N^{\mathbb{Q}_p}:\mathbb{Q}_p\to\{0\}\cup\mathbb{Q}_p^\times/(1+Np\mathbb{Z}_p)$ for the quotient map, sending 0 to 0. Write $\sigma_N^{\mathbb{R}}:\mathbb{R}\to\{+,0,-\}$ for the usual sign function on $\mathbb{R}$ sending 0 to 0. The maps $\sigma_N^{\mathbb{R}}$ and $\sigma_N^{\mathbb{Q}_p}$ play a similar role for preparing semi-algebraic functions and decomposing semi-algebraic sets; see Propositions 1.2 and 1.3 below. Note that the quotient $\mathbb{Q}_p^\times/(1+Np\mathbb{Z}_p)$ is isomorphic to $\mathbb{Z}\times(\mathbb{Z}_p/Np\mathbb{Z}_p)^\times$ by considering the map $\overline{\mathrm{ac}}_N(u\cdot p^\ell):=u \bmod Np\mathbb{Z}_p$ for $u$ in $\mathbb{Z}_p^\times$ and $\ell\in\mathbb{Z}$.

The following results are basic one-dimensional variants of cell decomposition and preparation. They capture part of their flavor and lie close to the definitions of both o-minimality and Hensel minimality, as we will see further below.

Proposition 1.2 (basic cell decomposition in dimension one). *Let $K$ be either $\mathbb{R}$ or $\mathbb{Q}_p$. Let $X\subset K$ and $f:X\to K$ be semi-algebraic and $N>0$ be an integer. Then there are a finite set $C\subset K$ and an integer $M>0$ such that $\sigma_N^K(f(x))$ depends only on the tuple $\sigma_M^K(x-c)_{c\in C}$ for $x$ in $K$.*

If $K$ is $\mathbb{R}$, it follows easily from Proposition 1.2 that any semi-algebraic set $X\subset K$ is a finite union of open intervals and singletons, which are called (unary) cells. Full cell decomposition is then a higher-dimensional version of Proposition 1.2, where the finite tuple $C$ moves nicely in the fibers of the family. We refer the reader to [61] for a general account on cell decomposition in the real, o-minimal case and to [49], [53], [54] for the initial work on $p$-adic cell decomposition.

In fact, for applications towards real and $p$-adic integration, it is important to see more finely than in Proposition 1.2 how $f(x)$ itself can be written in terms of $x-c$, where $c\in C$ is taken with minimal distance to $x$. This goes under the name preparation and forms part of Denef's $p$-adic cell decomposition in [53], [54].

To continue with the simple one-dimensional set-up, let us state a basic one-dimensional variant of such finer preparation.

Proposition 1.3 (basic preparation in dimension one). *Let $K$ be either $\mathbb{R}$ or $\mathbb{Q}_p$. Let $X\subset K$ and $f:X\to K$ be semi-algebraic. Then there are a finite set $C\subset K$, an integer $M>0$, and for each $c\in C$ and each $\lambda\in K$ a semi-algebraic function $u$ which is a unit, integers $r>0$ and $s$, and a constant $a\in K$ such that, for all $x\in K$ that lie closest to $c$ (among the elements of $C$) and which satisfy that $x-c$ lies in $\lambda\cdot P_M^K$, one has*

$$|f(x)|^r=|a|\cdot|u(x)|\cdot|x-c|^s,$$

*where $|\cdot|$ denotes the norm on $K$. Here, by a unit on a subset of $K$ we mean a function whose range is contained in a compact subset of $K^\times$.*

Much more can be said about preparation, also in families, and much more can be understood and assumed about the units $u$ in order to make it clear and simple to calculate an integral of $|f|$ over $X$ against the Haar measure (also in families), but we will just refer the reader to [96], [86], [51], [90], [45] for their study and initial uses concerning real integration and to [53] for $p$-adic integration.

**1.1** Before turning to more technical matters, let us give some further context related to Igusa's and Denef's work [78], [53], [56].

Consider a polynomial $f$ in $n$ variables with coefficients in $\mathbb{Z}$ and, for any integer $m>0$ and any prime number $p$, the finite exponential sum

$$S_f(p^m):=\frac{1}{p^{mn}}\sum_{\substack{x_j=1,\dots,p^m\\ j=1,\dots,n}}\exp\left(2\pi i\frac{f(x)}{p^m}\right).$$

Note that the sum $S_f(p^m)$ is normalized (it is divided by the number of summands), so that its complex modulus is at most 1. Igusa put forward a series of conjectures about the sums $S_f(p^m)$ that had a similar ambition and flavor as Weil's conjectures; see [78], [56].

Let us phrase a rationality result that is close in spirit to the rationality results of Igusa and Denef (see [78], [53]), but which took until 2005 [40], [39], [43] to be understood in greater generality, relying in particular on general variants of Propositions 1.2 and 1.3.

Theorem 1.4. *Let $f$ be a polynomial over $\mathbb{Z}$. Then*

$$\sum_{m>0}S_f(p^m)T^m$$

*is a rational function in $T$ over $\mathbb{C}$. That is, it lies in $\mathbb{C}(T)$.*

In fact, rationality results like Theorem 1.4 are now more finely understood by a general study of parameter integrals of functions of certain classes as we will see further below and which is also enabled by variants of cell decomposition and preparation.

Let us end the introduction with a result and a question about upper bounds for these exponential sums.

THEOREM 1.5 (Katz [82], extending Deligne [52]). *Let $f$ be a polynomial over $\mathbb{Z}$ in $n$ variables, and suppose that the top degree homogeneous part $f_d$ of $f$ is of degree $d>1$ and singular in (affine) dimension $s$. (Namely, the locus of* $\operatorname{grad} f_d = 0$ *has Zariski dimension $s$ in $\mathbb{C}^n$.) Then there is $c>0$ such that for any prime number $p$ one has*

$$|S_f(p)| \leq cp^{-\frac{n-s}{2}}.$$

CONJECTURE 1.6 (Igusa [78], variant from [47]). *Let notation be as in Theorem 1.5. Then there is $c>0$ such that for any prime number $p$ and any integer $m>0$ one has*

$$|S_f(p^m)| \leq cp^{-m\cdot\frac{n-s}{d}} \cdot m^n.$$

The independence of $c$ on $p$ and $m$ is of course key. Comparing Theorem 1.5 with this conjecture, one sees that the decay rate of $\frac{n-s}{2}$ for the finite field exponential sums conjecturally becomes $\frac{n-s}{d}$ in the $p$-adic case. This conjecture is open even for any fixed $p$, and even in the case $s=0$. Nguyen [94] recently found uniform upper bounds with decay rate $(n-s)/2(d-1)$ instead of $(n-s)/d$. The case of Conjecture 1.6 with $\frac{n-s}{d} \leq 1$ is covered by [46], [47], and the case with $n \leq 4$ is covered in [47]. The factor $m^n$ is of a logarithmic nature (in $p^m$) and is may not be necessary (see [47, sect. 1.4]).

Note that uniform bounds for the sums $S_f(p^m)$ play a role in Birch's work [12] on local global principles for polynomials $f$ of very low degree $d$ compared to $n-s$, namely for estimating the major arcs in the circle method. Upper bounds as in Conjecture 1.6 could play a role in relaxing Birch's relation from [12] up to $d<(n-s)/2$, albeit that not just the major arcs pose a challenge, but even more so do the minor arcs.

**2 From cell decomposition to o-minimality and Hensel minimality.** The notion of Hensel minimality (and its instance of 1-h-minimality) is an analogue of o-minimality for the non-archimedean setting and has been recently developed in [35] and [36]. These notions of 1-h- and of o-minimality are built upon insights coming from cell decomposition results (and their proofs and usages), and their definitions boil down to conditions on unary sets, as will soon be recalled.

Ideas of cell decompositions came to the foreground for real semi-algebraic sets via the work of Collins [50] (in his study of real cylindrical decompositions and related efficient algorithms) and by Cohen [49] and Denef [53], [54] in the $p$-adic case. They then evolved towards o-minimality in the real case [101], [61], towards Hensel minimality in the non-archimedean case [35], [36], and towards applications for integration (real, as well as $p$-adic and motivic) and for rational points, as we will see in the next sections. In this section, we recall and motivate the definitions of o- and Hensel minimality and some properties.

Let $\mathbb{R}$ be equipped with a structure coming from a language $\mathcal{L}$ containing the language of rings. Fixing such a structure is equivalent to fixing a collection of subsets of $\mathbb{R}^n$, for each $n$, which is stable under Boolean combinations and coordinate projections and which contains all zero loci of polynomials with rational coefficients (see [61] or [89] for further details). The sets appearing in this collection are called the definable sets, and a function between definable sets is called a definable function if its graph is a definable set. Let $\sigma^{\mathbb{R}} := \sigma_1^{\mathbb{R}}$ be as above. The following slightly altered definition of o-minimality resembles very closely the one-dimensional cell decomposition from Proposition 1.2 in the case of real semi-algebraic sets (which form an example of an o-minimal structure).

DEFINITION 2.1. *The $\mathcal{L}$-structure on $\mathbb{R}$ is called o-minimal if for every definable set $X \subset \mathbb{R}$ there is a finite set $C \subset \mathbb{R}$ such that for any $x \in \mathbb{R}$, "belonging to $X$" is determined by the tuple $(\sigma^{\mathbb{R}}(x-c))_{c\in C}$.*

Hensel minimality [35] and [36] is also inspired by Proposition 1.2 (in the $p$-adic case) but is more delicate for a couple of reasons. In fact, Hensel minimality works generally for a henselian valued field of characteristic zero and not just for $\mathbb{Q}_p$. Let us give a first very direct and basic analogue to Definition 2.1 for valued fields, in the equicharacteristic zero case for simplicity. Let $K$ be a henselian valued field of equicharacteristic zero, equipped with a structure coming from a language $\mathcal{L}$ containing the language of valued fields (that is, the language of rings together with a predicate for the valuation ring inside the valued field).

Definition 2.2 (almost-h-minimality). *To call $K$ almost-h-minimal, do as for o-minimality, with $\sigma^{\mathbb{R}}$ replaced by the quotient map*

$$\sigma^K : K \to K/(1+\mathcal{M}_K),$$

*where $\mathcal{M}_K$ is the maximal ideal of the valuation ring $\mathcal{O}_K$ of $K$ and where the quotient of semigroups is just the union of $\{0\}$ with $K^\times/(1+\mathcal{M}_K)$. In detail, call $K$ almost-h-minimal if for every definable set $X \subset K$ there is a finite set $C \subset K$ such that for any $x \in K$, "belonging to $X$" is determined by the tuple $(\sigma^K(x-c))_{c\in C}$.*

Often, the map $\sigma^K$ from Definition 2.2 is denoted by rv and the quotient $K/(1+\mathcal{M}_K)$ by RV (this notation indicates that it combines information from the residue field and the value group). It is sometimes called the leading term structure of $K$, it was first used in this context in [6], [7], and it is related to hyperfields.

Let us describe the subtleties involved in order to get a deeper analogue of Definition 2.1 instead of Definition 2.2. In Definition 2.1, some aspects are automatic which are not automatic in Definition 2.2. For example, in the real case, the finite set $C$ can be taken to be definable over the same set of parameters that are used to define $X$ and therefore works automatically very well in families. To put this in a slightly more technical but essentially equivalent way, an analogous formulation with a finite set $C$ as in Definition 2.1 automatically follows inside any elementary equivalent structure as soon as it holds inside $\mathbb{R}$. Second, in valued fields one needs to look somewhat deeper than just up to the level of $K/(1+\mathcal{M}_K)$ with the map $\sigma^K$ and even deeper than the above maps $\sigma_M^K$ used in Proposition 1.2, something which is not really relevant for the sign map on $\mathbb{R}$.

Let us now for the sake of precision give the definition of 1-h-minimality in the equicharacteristic zero case, which is an important form of Hensel minimality (see [35], [36], [108] for variants).

Write $\Gamma_K$ for the (multiplicatively written) value group of $K$, and for any $\lambda \le 1$ in $\Gamma_K$ write $\mathrm{rv}_\lambda$ for the quotient map $K \to \mathrm{RV}_\lambda$, with $\mathrm{RV}_\lambda$ being $K/(1+\lambda'\mathcal{M}_K)$, where $\lambda'$ is any element of $\mathcal{O}_K$ with valuation $\lambda$. Note that $\mathrm{RV}_1$ is just RV as introduced above.

Definition 2.3 (1-h-minimality [35]). *Let $K$ be a henselian valued field of equicharacteristic zero, equipped with an $\mathcal{L}$-structure as above.*

*Say that $K$ is 1-h-minimal if each elementary equivalent structure is weakly 1-h-minimal.*

*Say that $K$ is weakly 1-h-minimal if for any $\lambda \le 1$ in $\Gamma_K$, any $\xi \in \mathrm{RV}_\lambda$, and any $X$ which is definable using parameters from $\xi$, $A \subset K$, and RV, there is a finite $A$-definable set $C$ which prepares $X$ in the following way:*

*Whether $x \in K$ belongs to $X$ or not depends only on the tuple $\mathrm{rv}_\lambda(x-c)_{c\in C}$.*

Definition 2.3 makes sense for henselian valued fields of equicharacteristic zero only; for a generalization of this definition which includes the mixed characteristic case, we refer the reader to [35], [36].

That $C$ in Definition 2.3 is $A$-definable means that an $\mathcal{L}$-formula with parameters using elements from $A$ is allowed to describe $C$. The set $X$ in Definition 2.3 is described by an $\mathcal{L}$-formula using some elements from $A$, RV and possibly also $\xi$. (To use $\xi$ and elements of RV as parameters is slightly technical: it means that there is a formula describing $X$ in the expansion of $\mathcal{L}$ obtained from $\mathcal{L}$ by adding some predicates, namely for the inverse images in $K$ of $\xi$ and the occurring elements of RV; see [35].)

To indicate the strength and ambitions of Hensel minimality, we describe a partially known goal statement on Taylor approximation of any finite order and in general dimension.

Goal Statement 2.4 (full Taylor approximation (partially known)). *Let $K$ be a henselian valued field of equicharacteristic zero. Let $K$ be equipped with a structure which is 1-h-minimal, let $f : K^n \to K$ be definable, and fix $r > 0$. Then there exists a definable function $\chi : K^n \to \mathrm{RV}^m$ for some $m$ such that for every open fibre $F$ of $\chi$ we have that $f$ is $C^r$ on $F$, the norm of $f^{(i)}$ is constant on $F$ for all $i$ with $|i| \le r$, and for all $x, y \in F$ we have that*

$$|f(x) - T^{<r}_{f,y}(x)| \le \max_{|i|=r}(|f^{(i)}| \cdot |(x-y)^i|).$$

In Goal Statement 2.4, $T^{<r}_{f,y}$ denotes the Taylor polynomial of $f$ of order $< r$ at $y$ and the norm notation $|\cdot|$ is multiplicative notation for the valuation map on $K$. Further, $(x-y)^i$ uses multi-index notation, $f^{(i)}$ stands for the iterated derivative of $f$, and $|i| = i_1 + \cdots + i_n$.

Currently, the case with $n = 1$ is known, and more generally when $x$ and $y$ lie in a single box $B$ (namely a Cartesian product of balls) which is contained in a fiber of $\chi$; see [35, Theorems 3.2.2, 5.4.10, 5.6.1]. In the semi-algebraic and subanalytic frameworks (which are special cases of 1-h-minimality [38], [35]), Goal Statement 2.4 is more easy to obtain because of piecewise analyticity of definable functions. The general case of Goal Statement 2.4 is still a work in progress [37] at the moment of writing this text. Note that the finite tuple $(\mathrm{rv}(x-c))_{c\in C}$ as appearing in Definition 2.2 to prepare $X$ is replaced in higher dimensions by a definable map $\chi$ as in Goal Statement 2.4 to prepare $f$.

Historically, since the introduction of o-minimality, there has been a quest to find suitable analogue frameworks in the non-archimedean context, that is, for definable sets in henselian valued fields (of characteristic zero). A first studied structure expanding the $p$-adic semi-algebraic sets is formed by the $p$-adic subanalytic sets, introduced in [58] (following Hironaka [74], who introduced real subanalytic sets), for which cell decomposition was obtained in [91] and a finer preparation in [18]. One of the challenges of using an axiomatic framework to capture tame geometric properties is that analyticity can never be used as a tool since analyticity is not an axiomatic, first order property. Through the years there have been several ideas to capture tame non-archimedean geometry in axiomatic ways, but most of them seem either too restrictive or not powerful enough. Let us mention C-minimality [70], [88], P-minimality [71], v-minimality [75], and b-minimality [41], among other notions of tameness like NIP and dp-minimality [79]. Hensel minimality seems both versatile and powerful to be a full analogue of o-minimality for the non-archimedean context (that is, for studying definable sets and functions in characteristic zero henselian valued fields). Of course there are several other non-archimedean approaches to geometry, often more closely inspired by schemes and their morphisms than by o-minimality, as for example the spaces introduced by Berkovich [8] and [62] and a more model theoretic variant introduced by Hrushovski and Loeser in [76]. A model theoretic approach to tame geometry in positive characteristic henselian valued fields remains still very elusive and mysterious as far as we can see, although large positive characteristic can often be obtained from a uniform understanding of mixed characteristic as we will soon see below.

**3 From cell decomposition to *p*-adic, real, and motivic integration.** Cell decomposition and preparation results allow one to study parameter integrals of various kinds. Some of these results can be summarized by stating that certain classes of functions are stable under integration; namely, after integrating with respect to some of the variables one gets a function of the same class.

**3.1 *p*-adic and real integration.** Let us describe some systems of functions which are stable under integration and which are furthermore geometrically interesting in the $p$-adic and the real cases.

Definition 3.1. *For $X\subset\mathbb{Q}_p^n$ semi-algebraic, let $\mathcal{C}(X)$ be the $\mathbb{Q}$-algebra generated by the functions $|f|:X\to\mathbb{R}$ and $\log_p|g|:X\to\mathbb{R}$ for all semi-algebraic functions $f:X\to\mathbb{Q}_p$ and $g:X\to\mathbb{Q}_p^\times$.*

Here, the norm on $\mathbb{Q}_p$ sends $p$ to $p^{-1}$ and $\log_p$ is the logarithm with base $p$. Note that $-\log_p|x|$ is just the $p$-adic valuation of $x\in\mathbb{Q}_p^\times$. Denef was the first to study the algebras $\mathcal{C}(X)$ [55], [57], leading to the following more general result.

Theorem 3.2 (stability under integration [26]). *Let $f(x,y)$ be a function on $\mathbb{Q}_p^n\times\mathbb{Q}_p^m$ which is of class $\mathcal{C}$ (that is, $f$ lies in $\mathcal{C}(\mathbb{Q}_p^n\times\mathbb{Q}_p^m)$). Then there is $g$ in $\mathcal{C}(\mathbb{Q}_p^n)$ such that for any $x$ in $\mathbb{Q}_p^n$ one has*

$$g(x)=\int_{y\in\mathbb{Q}_p^m} f(x,y)|dy|$$

*whenever the right-hand side is finite (in the sense that the function $y\mapsto f(x,y)$ is integrable over $\mathbb{Q}_p^m$ for the normalized Haar measure $|dy|$ on $\mathbb{Q}_p^m$).*

Let us compare this Theorem 3.2 right away to the real case, where we switch to the (globally) subanalytic framework instead of the semi-algebraic one, a generalization initiated by Hironaka [74]. Note that integrals of real semi-algebraic functions are much harder and intimately related to the study of periods [85], [80], [81].

For $X\subset\mathbb{R}^m$ a subanalytic set, let $\mathcal{C}(X)$ be the algebra of functions $X\to\mathbb{R}$ generated by all subanalytic functions $f:X\to\mathbb{R}$ and all functions of the form $\log|g|$ with $g:X\to\mathbb{R}^\times$ subanalytic.

Theorem 3.3 (stability under integration [51], extended in [45]). *Let $f(x,y)$ be a function on $\mathbb{R}^n\times\mathbb{R}^m$, which is of class $\mathcal{C}$. Then there is $g$ in $\mathcal{C}(\mathbb{R}^n)$ such that for any $x$ in $\mathbb{R}^n$ one has*

$$g(x)=\int_{y\in\mathbb{R}^m} f(x,y)dy$$

*whenever the right-hand side is finite.*

In any case, functions of class $\mathcal{C}$ are real-valued, and they are obviously not stable under Fourier transforms. However, one can enrich them slightly so that they also contain functions of the form $x \mapsto \exp(2\pi i f(x))$ and so that they become stable under Fourier transform as well. This is done by forming the algebras of class $\mathcal{C}^{\exp}$ starting from $\mathcal{C}$. The class $\mathcal{C}^{\exp}$ is used in [1] and [2] to answer some questions of [3] about holonomicity of distributions on local fields (see also some further open questions at the end of [2]). In the real case. this is slightly more subtle than in the $p$-adic case, where the functions of class $\mathcal{C}^{\exp}$ can be easily defined as follows.

Choose a nontrivial additive character $\psi : \mathbb{Q}_p \to \mathbb{C}^\times$. (Namely, $\psi$ is a nontrivial continuous group homomorphism.)

Definition 3.4. *For $X \subset \mathbb{Q}_p^n$ semi-algebraic, let $\mathcal{C}^{\exp}(X)$ be the algebra generated by*
- *all functions in $\mathcal{C}(X)$ and*
- *all functions of the form $\psi(h)$ for semi-algebraic $h : X \to \mathbb{Q}_p$.*

Theorem 3.5 (stability under integration and Fourier [43], extended in [26]). *The class $\mathcal{C}^{\exp}$ is stable under integration and thus under Fourier transformation. In detail, let $f(x,y)$ be a function on $\mathbb{Q}_p^n \times \mathbb{Q}_p^m$ which is of class $\mathcal{C}^{\exp}$. Then there is $g$ in $\mathcal{C}^{\exp}(\mathbb{Q}_p^n)$ such that for each $x$ in $\mathbb{Q}_p^n$*

$$g(x) = \int_{y \in \mathbb{Q}_p^m} f(x,y)|dy|$$

*whenever the right-hand side is finite (in the sense that the function $y \mapsto f(x,y)$ is integrable over $\mathbb{Q}_p^m$).*

Note that integrability for the integral in Theorem 3.5 is rather subtle and harder to control than in Theorems 3.2 and 3.3 because of the oscillation coming from $\psi$. Other subtle integrability questions for a subclass of $\mathcal{C}$ are treated in [48]; these questions from [48] remain completely mysterious in the oscillatory context. The stability given by Theorem 3.5 is important for applications to the Langlands program, as we will see in section 4. In the real case, the $\mathcal{C}^{\exp}$-class is introduced in [21], where results like stability under integration are obtained.

Let us give one more result which depends very delicately on integrability because of oscillation. It gives a close link between zero loci and loci of integrability for functions of $\mathcal{C}^{\exp}$-class. Such results exist in several other settings (for example, uniformly in $p$ as sketched in the next subsection), but in the real setting of [21] this is still open for the $\mathcal{C}^{\exp}$ class, except in some low dimensions.

Theorem 3.6 (correspondence of zero loci with loci of integrability [26]). *Let $f(x,y)$ be a function on $\mathbb{Q}_p^n \times \mathbb{Q}_p^m$ which is of class $\mathcal{C}^{\exp}$. Then there is $g$ in $\mathcal{C}^{\exp}(\mathbb{Q}_p^n)$ such that the zero locus of $g$ inside $\mathbb{Q}_p^n$ equals the following locus of integrability:*

$$\{x \in \mathbb{Q}_p^n \mid \text{ the function } y \mapsto f(x,y) \text{ is } L^1 \text{ on } \mathbb{Q}_p^m\},$$

*and vice versa (for each $g$, there exists $f$).*

Theorem 3.6 can be seen as relating analysis and geometry. Zero loci of functions of $\mathcal{C}$-class or of $\mathcal{C}^{\exp}$-class are more general than semi-algebraic sets, but they are related to their study and have some stratifications and dimensions; see [34], [1]. Most of the above $p$-adic concepts and results make sense uniformly in $p$ by more recent (uniform) analogues of Macintyre's quantifier elimination (mutatis mutandis). They hold, moreover, uniformly over all finite field extensions of $\mathbb{Q}_p$ for all $p$ and even often over all fields $\mathbb{F}_q((t))$ for $q$ any power of a *large enough* prime. This can be seen in at least two ways: one can boil all aspects down to a use of Gödel's completeness theorem [64], [65], [89] so that only finitely many positive characteristics can be excluded, or one ensures all the way through that results work for a large enough positive characteristic. For both ways, it is important to have building blocks that are definable in a uniform way. We will give a sketch of this in the next subsection.

**3.2 Uniform $p$-adic and motivic integration, a brief sketch.** As already mentioned, most of the above results about $p$-adic integration can be made uniform in the local field and some hold even at a deeper, motivic level. This needs uniform notions for semi-algebraic (or definable) sets and for the rings of functions of classes $\mathcal{C}$ and $\mathcal{C}^{\exp}$ on them. In the uniform $p$-adic and the motivic cases, the definitions of the classes $\mathcal{C}$ and $\mathcal{C}^{\exp}$ are not as short and simple as for a fixed $p$-adic field in Definitions 3.1 and 3.4. Let us highlight some of the important ingredients, first for the uniform $p$-adic setting, following [32] and [29]. Another approach to uniform $p$-adic integration uses embedded resolutions with good reduction (instead of cell decomposition; see [56]), but that approach does not work so well for family settings, Fubini, and Fourier transforms.

Let us give some ingredients and ideas for the adaption to the uniform $p$-adic setting, without going in full detail. By a non-archimedean local field we mean a finite field extension of $\mathbb{Q}_p$ or of $\mathbb{F}_p((t))$ for some prime number $p$. Let Loc be the collection of pairs $(K,\varpi)$ where $K$ is a non-archimedean local field and $\varpi$ is a uniformizer of $\mathcal{O}_K$. The semi-algebraic sets as defined in the introduction need to be generalized slightly so that they work well uniformly over the local fields. To this end, the generalized language of Denef-Pas $\mathcal{L}_{\rm DP}$ from [32] is used (generalizing the setups from [97], [98] for uniform semi-algebraic $p$-adic sets). This leads to a notion of definable sets $X$ given by $\mathcal{L}_{\rm DP}$-formulas. Such a definable set $X$ is in essence a collection of sets

$$X=(X_K)_{K\in \mathrm{Loc}}$$

with $X_K$ now being a subset of

$$K^n\times\mathbb{Z}^m\times(\mathcal{O}_K/N\mathcal{M}_K)^r$$

for some nonnegative integers $n,m,r$ and some integer $N>0$. (The uniformizers of the pairs $(K,\varpi)$ in Loc are often implicit and left out from the notation.) An important aspect of the sets $X_K$ is that the $n,m,r,N$ do not change with $K$; indeed, $X$ is required to be described by an $\mathcal{L}_{\rm DP}$-formula.

For the class $\mathcal{C}$, there is an important new generator coming from counting the number of elements inside definable subsets of the finite rings of the form $(\mathcal{O}_K/N\mathcal{M}_K)^r$, namely mapping $x$ in $X_K$ to the number

$$\#Y_{K,x} \tag{3.1}$$

for some definable set $Y$ with $Y_K\subset X_K\times(\mathcal{O}_K/N'\mathcal{M}_K)^{r'}$ and $Y_{K,x}$ the fiber above $x$; note that the counting thus happens inside the finite set $(\mathcal{O}_K/N'\mathcal{M}_K)^{r'}$.

For the class $\mathcal{C}^{\rm exp}$, there are two further key changes: the additive character $\psi$ on $K$ is allowed to vary over all characters with fixed conductor (namely, constant on $\mathcal{M}_K$ and nonconstant on $\mathcal{O}_K$) and the above counting function from (3.1) is replaced by the following form of finite exponential sum function, sending $x$ in $X_K$ to

$$\sum_{y\in Y_{K,x}}\psi(h(x,y)), \tag{3.2}$$

where $h=(h_K)_K$ is a collection of functions $h_K:Y_K\to K/\mathcal{M}_K$ coming from a definable function in the sense that there is a definable subset $Z$ with $Z_K\subset Y_K\times K$ such that $Z_K$ projects to the graph of $h_K$ for each $K$.

With these key changes in mind (and omitting some of the details of the definitions), we can now state the main results for the classes $\mathcal{C}$ and $\mathcal{C}^{\rm exp}$ in the uniform case. First, we state a uniform version of Theorems 3.2 and 3.5.

THEOREM 3.7 (stability under integration [26], [32]). *Let $f=(f_K)_K$ be of class $\mathcal{C}$, resp. $\mathcal{C}^{\rm exp}$, with domain $X\times Y$ for definable sets $X=(X_K)_K$, $Y=(Y_K)_K$. Then there are $g$ of class $\mathcal{C}$, resp. $\mathcal{C}^{\rm exp}$, and an integer $M>0$ such that, as soon as the characteristic of $K$ is zero or at least $M$, one has for all $x$ in $X_K$ that*

$$g_K(x)=\int_{y\in Y_K}f_K(x,y)|dy|$$

*holds whenever the right-hand side is finite, with $|dy|$ the product measure of the Haar measure and the counting measure.*

As one sees in Theorem 3.7, to make sensible statements, one has to exclude small positive characteristics. This is even the case for rationality results as in Theorem 1.4, where the (small) positive characteristic case is still a wide open question; let us refer the reader to [77], [93] for various generalizations of rationality results. The correspondence of zero loci with loci of integrability (similar to Theorem 3.6) also holds uniformly in $K$ (again after excluding some positive characteristics for the local field); see [26], [29].

In the next section, we will see that there are close relations between integrals and integrability for local fields with the same residue field, as soon as the residue field characteristic is large enough. This leads to transfer principles that have been applied in the Langlands program.

For genuine motivic integration and stability under integration, Fubini, and Fourier as developed in [42], [43], a key point is to generalize the generators (loosely) described in (3.1) and (3.2) further by replacing the counting (resp. the exponential sum) functions by the isomorphism class of $Y$ in a certain Grothendieck ring of definable sets (resp. of pairs $(Y, h)$), seen relatively, namely as a family over $X$. This gives a form of stability under Fourier transforms at a deeper level than just uniformly in $K$. These frameworks have been extended to also include Mellin transforms by providing some abstract avatars at the level of Grothendieck rings for multiplicative characters in [44], and in the real case the Mellin transform is added in [22], [23]. Note that at this deeper motivic level it is harder to define the notions of integrability and of loci, some aspects of which still need to be developed, with a recent step made in [33]. For a variant of motivic integration without Fubini and Fourier, we refer the reader to the book [17], which follows more closely the frameworks from [84], [59], [60]. For a geometric variant of motivic integration working inside algebraically closed valued fields also based on cell decomposition and with Fubini and Fourier, we refer the reader to [75], and to [106] for a generalization of this setup to the Hensel minimal setting (inside not necessarily algebraically closed valued fields).

**4 Transfer principles for uniform $p$-adic integration and the fundamental lemma.** With the setup from subsection 3.2, let us right away state two key transfer principles (further transfer principles can be found in [28], [29]).

THEOREM 4.1 (transfer principle for equalities [43]). *Let $F$ and $G$ be functions in $n$, resp. in $m$, valued field variables and of class $\mathcal{C}^{\exp}$ (as in subsection 3.2). Then there exists a lower bound $M$ such that for any local field $K$ (in Loc) with residue field characteristic at least $M$ the following holds.*

*Whether or not*

$$\int_{K^n} F_K = \int_{K^m} G_K \text{ holds for each } \psi$$

*only depends on the residue field of $K$.*

Note that the additive characters $\psi$ are quite different over $\mathbb{Q}_p$ than over $\mathbb{F}_p((t))$ and have a totally different image in the complex unit circle, which makes it hard at first to change the characteristic. An orthogonality argument comes into play to show Theorem 4.1; see [43], [44].

Without going into any details, let us state its first application to the Langlands program.

THEOREM 4.2 ([30]). *The fundamental lemma of the Langlands program follows in characteristic zero (from Ngô's and Hales' works).*

Indeed, the fundamental lemma is about an equality of $p$-adic integrals that turn out to be of class $\mathcal{C}^{\exp}$, and Ngô famously obtained these equalities in positive characteristic [92]. By Theorem 4.1 the equalities then follow in characteristic zero with large residue field characteristic, and by [69] they then also follow for the remaining small residue field characteristics. This represents an alternative strategy to work by Waldspurger [110] [111], with the upshot that it also works for variants of the fundamental lemma (a weighted variant, a variant by Jacquet and Ye, and further variants in [114] and [16]; see also [66]). Some twisted variants of the fundamental lemma remain to be put in the framework of the class $\mathcal{C}^{\exp}$ up to date as far as we can see.

THEOREM 4.3 (transfer principle for integrability [26]). *Let $F$ be a function of class $\mathcal{C}^{\exp}$ (as in subsection 3.2). Then there exists a lower bound $M$ such that for any local field $K$ (in Loc) with residue field characteristic at least $M$ the following holds. Whether or not $F_K$ is integrable (for each $\psi$) only depends on the residue field of $K$.*

In [26], Theorem 4.3 is reduced to Theorem 4.2 by showing a correspondence between loci of integrability and zero loci for functions of $\mathcal{C}^{\exp}$-class, a uniform version of Theorem 3.6. As an application of Theorem 4.3, we obtained that the Harish–Chandra character of an irreducible admissible representation of $G(\mathbb{F}_q((t)))$ with $G$ a connected, reductive group is everywhere locally integrable for $p$ large enough (given the root datum of $G$) and $q$ a power of $p$ (under a slight extra condition of existence of the mock exponential map) [27]. This result uses the full power of $\mathcal{C}^{\exp}$ (in particular its stability under Fourier transforms), as the character was expressed (near tame semisimple elements) by its local character expansion as a finite linear combination of Fourier transforms of orbital integrals.

We refer the reader to [105] for applications of other transfer principles, for example for bounds of (the complex modulus of) functions of $\mathcal{C}^{\exp}$ class.

Let us end this section by mentioning that a genuine motivic version of the fundamental lemma has been obtained (at the level of Grothendieck rings) by Forey, Loeser, and Wyss [63], thereby completing a research program initiated by Hales [68].

**5 Rational points, Diophantine results, and a question.** In this section, we give some Diophantine results in Hensel minimal structures, some of which are still work in progress [37] at the moment of writing this text. We will end the section with a recent open question around tight dimension growth bounds for rational points on projective varieties.

Let us first recall the famous Pila–Wilkie counting theorem on rational points on definable sets in o-minimal structures.

THEOREM 5.1 ([100]). *Let $X \subset \mathbb{R}^n$ be definable in an o-minimal structure. Then for every $\varepsilon > 0$ there exists $c = c_\varepsilon$ such that for every $H \geq 1$ one has*

$$\# X^{\mathrm{trans}}(\mathbb{Q}, H) < cH^\varepsilon.$$

Here, $X^{\mathrm{trans}}(\mathbb{Q}, H)$ is the set of points $(x_1, \ldots, x_n)$ in $\mathbb{Q}^n$ which lie on $X^{\mathrm{trans}}$ and with each $x_i$ of height at most $H$, and $X^{\mathrm{trans}}$ is the set $X \setminus X^{\mathrm{alg}}$ where $X^{\mathrm{alg}}$ is the algebraic part of $X$. Recall that $X^{\mathrm{alg}}$ is the set of $x \in X$ for which there exists a semi-algebraic curve $C$ which contains $x$, which lies inside $X$, and which is of constant local dimension 1.

The Pila–Wilkie theorem (Theorem 5.1) plays an important role in many arithmetic applications, culminating quite recently in the solution of the André–Oort conjecture in [99]. A key role in Theorem 5.1 is played by the general theory of o-minimal structures, which recently found powerful applications in Hodge theory and the study of Hodge loci [5], [4], [83] as well.

Let us right away state the goal statement of a full non-archimedean analogue of Theorem 5.1 for Hensel minimal structures instead of o-minimal structures. This uses the mentioned generalization of Definition 2.3, which includes mixed characteristic from [35], [36].

GOAL STATEMENT 5.2 (Pila–Wilkie-type bounds for 1-h-minimal structures (partially known)). *Let $K$ be a finite field extension of $\mathbb{Q}_p$ for some prime number $p$. Let $X \subset K^n$ be definable in a 1-h-minimal structure on $K$. Then for every $\varepsilon > 0$ there exists $c = c_\varepsilon$ such that for every $H \geq 1$ one has*

$$\# X^{\mathrm{trans}}(\mathbb{Q}, H) < cH^\varepsilon.$$

As Goal Statement 5.2 is still work in progress [37] at the moment of writing this text, let us make precise which special cases are already in the literature. The case of Goal Statement 5.2 where $X$ has dimension one is fully obtained in [36] and is further generalized in [15], [109]. Goal Statement 5.2 is also known in the special case that $X$ is a subanalytic set of any dimension [20] and holds even uniformly in $p$, implying also the case of large positive characteristic [24] (see [11] for a study of rational points on zero loci of analytic functions that works in any positive characteristic). The missing ingredient in [36] to treat the general case of Goal Statement 5.2 is in particular a good understanding of Taylor approximation in general dimension (as in Goal Statement 2.4) and parametrizations with good Taylor approximation in general dimension. The latter is a non-archimedean analogue of Yomdin–Gromov parametrizations [67, 113, 112] and of the o-minimal parametrizations from [100]. In the present axiomatic setting of Hensel minimality, another challenge arises, namely that of bounding derivatives given that a function itself is nicely bounded. In the real setting, this comes more or less directly from the fundamental theorem of calculus, which allows one to relate the sizes of a monotone function and its derivative. In the analytic settings of [20, 24, 95], it comes from the relation between the Gauss norm (of a converging power series) and the supremum norm. In the one-dimensional case from [36], it comes from a Jacobian property which holds in Hensel minimal structures [35], [36]. Another difficulty is of course the total disconnectedness by which it is a challenge to go from local to piecewise, which was overcome in a first step for Lipschitz continuity in [19], [31].

**5.1 Towards tight projective dimension growth.** Let us end this paper with a discussion about rational points of bounded height on projective varieties. Projective dimension growth results were an open question raised by Serre and Heath-Brown in [104] and in [72], and it was answered mostly by Salberger [102], strengthening [14]. These results have subsequently found many generalizations and variants, of which we state the recent one with a quadratic dependence on the degree from [9], together with worst case lower bounds from [25].

For $X \subset \mathbb{P}^n_{\mathbb{Q}}$ and any $H \geq 1$, write $X(\mathbb{Q}, H)$ for the set of rational points $x$ on $X$ that can be written with integer homogeneous coordinates $(x_0, \ldots, x_n)$ where the integers $x_i$ satisfy $|x_i| \leq H$ for all $i$.

Theorem 5.3 (projective dimension growth [9], lower bounds from [25]). *Let $n > 1$ be given. Then there are constants $c$ and $\kappa$ such that for any integral hypersurface $X$ of degree $d \geq 4$ in $\mathbb{P}^n_{\mathbb{Q}}$ and any $H > 2$ one has*

$$\#X(\mathbb{Q}, H) \leq cd^2 H^{\dim X} (\log H)^{\kappa}. \tag{5.1}$$

*Furthermore, there is a constant $c' > 0$ such that for arbitrarily large values of $d$ and $H$ there exists an integral hypersurface $X$ of degree $d$ in $\mathbb{P}^n_{\mathbb{Q}}$ such that*

$$c' d^{2-2/n} H^{\dim X} \leq \#X(\mathbb{Q}, H). \tag{5.2}$$

As one can see, the worst case lower bound from (5.2) is still some distance away from the general upper bounds from (5.1), as far as the exponent of $d$ is concerned. The author and Glazer [25] raise the question of whether the factor $d^2$ in (5.1) can be improved to match the $d^{2-2/n}$ from (5.2), which is open even in the curve case (with $n = 2$), where the conjectured upper bound would become $cdH$ (omitting logarithmic factors for a moment). Note that in the curve case general upper bounds of the form $cd^2H^{2/d}(\log H)^{\kappa}$ are known for all $d \geq 1$ by [10] (this improves bounds from [13], [73]) and are of a flavor somewhat different from the conjectural bounds of size $cdH$, which can be compared with the basic Schwartz–Zippel bounds of the form $cdH^2$ in this projective case with $n = 2$. (Note that Schwartz–Zippel bounds also hold in the reducible case.) The proof of Theorem 5.3 does not rely on cell decomposition techniques, but rather on variants of the $p$-adic determinant method from [102], initiated in [73], and intersection theory. By the simplicity and generality of these suggested bounds, we believe that this question represents an interesting challenge for future research.

**Acknowledgments.** The author would like to express his gratitude to Jan Denef and François Loeser for so many stimulating discussions and encouragement throughout the years. The author would equally like to thank all of his co-authors for such nice and stimulating collaborations. The author was partially supported by KU Leuven IF C16/23/010 and by the support of the CDP C2EMPI, as well as the French State under the France-2030 programme, the University of Lille, the Initiative of Excellence of the University of Lille, the European Metropolis of Lille, and the R-CDP-24-004-C2EMPI project.